\documentclass[a4paper,10pt]{article}

\usepackage{amssymb, amsmath, amsthm, eucal, amscd}
\usepackage{pictexwd,dcpic}

\newtheorem{teorema}{Theorem}[section]
\newtheorem{definicao}[teorema]{Definition}
\newtheorem{lema}[teorema]{Lemma}
\newtheorem{proposicao}[teorema]{Proposition}
\newtheorem{corolario}[teorema]{Corolary}

\newenvironment{prova}{\setlength{\parindent}{0pt}\textbf{Proof.}}{\hspace{\stretch{1}}$\Box$}

\title{\textbf{On a paraconsistentization functor in the category of consequence structures}}

\author{Edelcio G. de Souza$^{1}$ \\ 
Alexandre Costa-Leite$^{2}$ \\
Diogo H.B. Dias$^{1}$ \\
\\
\footnotesize $^{1}$ University of Sao Paulo (BR) \ $^{2}$ University of Brasilia (BR)}

\date{ }

\begin{document}

\maketitle

\begin{abstract}
This paper is an attempt to solve the following problem: given a logic, how to turn it into
a paraconsistent one? In other words, given a logic in which \emph{ex falso quodlibet} holds,
how to convert it into a logic not satisfying this principle? We use a framework
provided by category theory in order to define a category of consequence structures. 
Then, we propose a functor 
to transform a logic not able to deal with contradictions into a paraconsistent one. 
Moreover, we study the case of paraconsistentization of propositional classical logic.
\end{abstract}


\section{Introduction}

Paraconsistent logics were initially studied by S. Ja\'skowski in \cite{jaskowski} 
and N. da Costa in \cite{dacosta}.
Since then, much research has changed the scenario in this field. 
Lately, we have many (paraconsistent) logics but also many distinct 
approaches to them (for a detailed survey on these inquiries, systems and methodologies,
see \cite{priest} and also \cite{dacosta-krause-bueno}). 
Some argued that the definition of a paraconsistent logic is not clear and sufficiently restrictive, 
because it is not able to determine univocally 
what a paraconsistent logic is (see \cite{beziau1}). Nevertheless, it is standard to accept these logics as those denying some form of the
principle known as \emph{ex falso quodlibet}, for which, given a formula schema $\varphi$ and its negation $\neg \varphi$, 
anything follows (see \cite{beziau}):

\begin{center}
$\varphi, \neg \varphi \vdash \psi$
\end{center}

This principle - also widely known as \emph{principle of explosion} - has a variety of formulations (just to consider a few examples):

\begin{center}
$\varphi \wedge \neg \varphi \vdash \psi$\\
$\vdash \varphi \rightarrow (\neg \varphi \rightarrow \psi)$\\
$\vdash (\varphi \wedge \neg \varphi) \rightarrow \psi$
\end{center}

These forms of \emph{ex falso} hold in classical logic. But, in most cases, in a paraconsistent
logic, they should fail.

The problem examined in this paper is: given a logic in which \emph{ex falso} holds, how can it be converted 
into a paraconsistent one? Or, is there a procedure for \emph{paraconsistentizing} a logic? Formally, given
a logic $L_1$ such that $\varphi, \neg \varphi \vdash_{L_1} \psi$, we have to find a procedure to turn
$L_1$  into a logic $L_2$ such that  $\varphi, \neg \varphi \nvdash_{L_2} \psi$. Further, if $L_1$ has a set 
of properties $C$, then, does a paraconsistent version of it, let's say, $L_2$,
preserve some of these properties?

The problem above combined with the idea of a \emph{paraconsistentization}
has been addressed by Costa-Leite in \cite{costaleite}. However, the author did not present a general and unified
method for producing paraconsistentization, but rather he proposes only particular ways of turning classical logic and modal logics into
paraconsistent logics. Previously, Beziau in \cite{beziau1} showed how to get some paraconsistent logics from modal 
logics using translations of logics, but this is not a method to convert any given
logic into a paraconsistent one. In the spirit of the goal of this
paper, combinations of paraconsistent with modal logics convert these into logics capable of tolerating inconsistencies.
This approach has been developed by people working with fibring logics, especially in the paper \cite{caleiro}.
Afterwards, Payette in \cite{payette} explored a way of generating inconsistent non-trivial logics from
consistent ones using a variation
of forcing. Differently, Caminada, Carnielli and Dunne studied in \cite{carnielli} a semi-stable semantics
able to find criteria for a given formal system to have paraconsistent characteristics. 

Our approach uses tools from category
theory and general abstract logic in order to provide a way to 
paraconsistentize any logic. The aim is to get a general, unified and abstract
perspective by means of a functor from the category of explosive logics 
(accepting \emph{ex falso}) to the category of 
non-explosive logics (\emph{rejecting ex falso}). The idea of using
category theory to reason about logic is very popular recently, 
but studies connecting paraconsistency and categories are rare, and we
could mention, for instance, the case of \cite{vasyukov} and also \cite{mortensen}.

In order to give an answer to the problem raised in \cite{costaleite}, we
start by defining a category of consequence structures (logics in a very
abstract sense) called $CON$, and we examine some remarkable features of it. 
Then, we define a functor using the category $CON$ and show some properties preserved
by this functor (the paraconsistentization functor). We proceed by defining paraconsistent
consequence structures and presenting some sufficient conditions to convert a given
logic into a paraconsistent one. By the end of this article, a particular paraconsistentization is presented.
We take the case of classical propositional logic and show how to paraconsistentize it.


\section{The category $CON$}

Alfred Tarski, in \cite{tarski} and \cite{tarski1}, came up with a definition of logical consequence (i.e. consequence operator) which allowed logicians to
reason at the abstract level, characterizing this notion by some conditions which are known as \emph{Tarskian}. 
In this paper, we do not impose any condition on the consequence operator $Cn$ (e.g., inclusion, idempotency, monotonicity, finiteness and so on). 
Tarskian structures are particular cases of consequence structures.\footnote{Logical structures without axioms or restrictions are proposed in \cite{beziau2}.}

A \emph{consequence structure} is a pair $(X,Cn)$ such that $X$ is a set and $Cn$ is an operation in $\wp(X)$, the power set of X:
$$
Cn: \wp(X) \rightarrow \wp(X).
$$

If $A$ is a subset of $X$, $Cn(A)$ is the set of $Cn$-\emph{consequences} of $A$ in the structure $(X,Cn)$. The set $X$ is called the \emph{domain} of the structure and $Cn$ is its \emph{consequence operator}.

We say that a subset $A \subseteq X$ is $Cn$-\emph{consistent} if $Cn(A) \neq X$; otherwise, $A$ is called $Cn$-\emph{inconsistent}. 

In what follows, if $f:X \rightarrow X'$ is a function, we use also $f: \wp(X) \rightarrow \wp(X')$ to denote the function $f$ extended to the power sets. If $A \subseteq X$, we have:
$$
f(A) = \{f(a): a \in A\}.
$$

Let $(X,Cn)$ and $(X',Cn')$ be two consequence structures. A \emph{homomorphism} $h$ from $(X,Cn)$ to $(X',Cn')$, denoted by $(X,Cn) \stackrel{h}{\longrightarrow} (X',Cn')$, is a function $$h: X \rightarrow X'$$ such that:

(i) $h$ is 1-1 (an injection);

(ii) $h$ preserves the consequence operator, that is, the following diagram is commutative:

$$
\begindc{\commdiag}[100]
\obj(0,0)[c]{$\wp(X')$}
\obj(1,0)[d]{$\wp(X')$}
\obj(0,1)[a]{$\wp(X)$}
\obj(1,1)[b]{$\wp(X)$}
\mor{c}{d}{$Cn'$}[-1,0]
\mor{a}{c}{$h$}[-1,0]
\mor{b}{d}{$h$}
\mor{a}{b}{$Cn$}
\enddc
$$

Thus, $h \circ Cn = Cn' \circ h$, that is,  $h(Cn(A)) = Cn'(h(A))$, for all $A \subseteq X$.

Since $h$ is an injection, we have that \emph{homomorphisms preserve consistent sets}. 

\begin{lema}
\label{l21}
If $(X,Cn) \stackrel{h}{\longrightarrow} (X',Cn')$ is a homomorphism and $A \subseteq X$ is $Cn$-consistent, then $h(A)$ is $Cn'$-consistent.
\end{lema}

\begin{prova}
Suppose that $A \subseteq X$ is $Cn$-consistent and $h(A)$ is $Cn'$-inconsistent. So, $Cn'(h(A)) = X'$. Since $h$ is a morphism, $h(Cn(A)) = X'$. But $h$ is an injection. Therefore, $Cn(A) = X$ (contradiction!).
\end{prova}

\medskip

Compositions of homomorphisms are homomorphisms and the identity function on $X$ is a homomorphism. Compositions satisfy associativity and identities satisfy the identity laws. Therefore, we have \emph{the category of consequence structures,} denoted by $CON$, whose $CON$-objects are consequence structures and $CON$-morphisms are homomorphisms.

We analyze some properties of $CON$.

\begin{proposicao}
\label{prop22}
$CON$ is neither (finitely) complete nor (finitely) co-complete.
\end{proposicao}

\begin{prova}
We prove that $CON$ does not have limits and co-limits for the empty diagram.

Suppose that $CON$ has terminal object $(X, Cn)$ and $card(X)=\kappa$. Consider $(Y, Cn')$ such that $card(Y)>\kappa$. Since there is no injection from $Y$ to $X$, there is no morphism in $CON$ from $(Y, Cn')$ to $(X, Cn)$. But this is in contradiction with the initial supposition. Therefore, $CON$ does not have terminal object.

In $SET$, the category of sets, $\emptyset$ is an initial object. If $CON$ would have initial object, its domain should be $\emptyset$. But, in this case, the consequence operator would be $Id_{\{\emptyset\}}$, that is an injection. On the other hand, a $CON$-morphism $h$ from $(\emptyset,Id_{\{\emptyset\}})$ to $(X,Cn)$ would be such that $h(\emptyset) = \emptyset$.

But, it is easy to see that, in this case, the following diagram

$$
\begindc{\commdiag}[100]
\obj(0,0)[c]{$\wp(X)$}
\obj(1,0)[d]{$\wp(X)$}
\obj(0,1)[a]{$\wp(\emptyset)$}
\obj(1,1)[b]{$\wp(\emptyset)$}
\mor{c}{d}{$Cn$}[-1,0]
\mor{a}{c}{$h$}[-1,0]
\mor{b}{d}{$h$}
\mor{a}{b}{$Id_{\{\emptyset\}}$}
\enddc
$$
commutes only if $Cn(\emptyset) = \emptyset$.
\end{prova}

\medskip
By the way, category $CON$ has no co-products:\footnote{We thank to the anonymous referee for the proof
reproduced here.} Assume any two consequence structures $\mathcal{C} = (C,Cn_{C})$ and $\mathcal{D} = (D,Cn_{D})$ as well as their co-product $\mathcal{C} + \mathcal{D}$. There are then injections $inl: \mathcal{C} \rightarrow \mathcal{C} + \mathcal{D}$ and $inr: \mathcal{D} \rightarrow \mathcal{C} + \mathcal{D}$ providing for any $f: \mathcal{C} \rightarrow \mathcal{X}$ and $g: \mathcal{D} \rightarrow \mathcal{X}$ an $[f,g]: \mathcal{C} + \mathcal{D} \rightarrow \mathcal{X}$ with $[f,g] \circ inl = f$ and $[f,g] \circ inr = g$. Consider that $f,f': \mathcal{C} \rightarrow \mathcal{X}$ and $g,g': \mathcal{D} \rightarrow \mathcal{X}$ are such that there are $c \in C$, $d \in D$ with $f(c) = g(d)$ and $f'(c) \neq g'(d)$. We have that $f(c) = g(d)$ guarantees that $([f,g] \circ inl)(c) = ([f,g] \circ inr)(d)$; and given that $[f,g]$ must be 1-1, this requires $inl(c) = inr(d)$. But $f'(c) \neq g'(d)$ ensures that $([f',g'] \circ inl)(c) \neq ([f',g'] \circ inr)(d)$; and since $[f',g']$ is a function, it follows that $inl(c) \neq inr(d)$ (contradiction!).

The task to determine whether $CON$ has some other universal constructions should still be explored in detail, and they are beyond the scope of the present paper. We intend also to study the case in which $CON$-morphisms are not injections but only preserve consistent sets.


\section{The functor $\mathbb{P}$}

Let us construct an endofunctor $\mathbb{P}$ on the category $CON$ that will be called \emph{paraconsistentization functor}.\footnote{An initial step towards the construction of this functor has been proposed in \cite{dias}. Here we examine this functor in detail, presenting a full characterization of its
structure and showing properties preserved by it.}

If $(X,Cn)$ is a consequence structure, we define a new operation $Cn_{\mathbb{P}}: \wp(X) \rightarrow \wp(X)$ such that, for all $A \subseteq X$:
$$
Cn_{\mathbb{P}}(A) := \bigcup\{Cn(A'): A' \subseteq A,\   Cn\text{-consistent} \}.
$$
In this way, we have that $x \in Cn_{\mathbb{P}}(A)$ if and only if there exists $A' \subseteq A$ $Cn$-consistent such that $x \in Cn(A')$.

Now, we define the action of $\mathbb{P}$ on $CON$:

\begin{itemize}
	\item For $CON$-objects $(X,Cn)$, $\mathbb{P}(X,Cn) = (X,Cn_{\mathbb{P}})$;
	\item For $CON$-morphisms $h$, $\mathbb{P}(h) = h$.
\end{itemize}

\begin{proposicao}
\label{p31}
$\mathbb{P}$ is an endofunctor in the category $CON$.
\end{proposicao}

\begin{prova}
Consider the following diagram:

$$
\begindc{\commdiag}[120]
\obj(0,0)[c]{$(X,Cn_{\mathbb{P}})$}
\obj(1,0)[d]{$(X',Cn'_{\mathbb{P}})$}
\obj(0,1)[a]{$(X,Cn)$}
\obj(1,1)[b]{$(X',Cn')$}
\mor{c}{d}{$\mathbb{P}(h) = h$}[-1,0]
\mor{a}{c}{$\mathbb{P}$}[-1,0]
\mor{b}{d}{$\mathbb{P}$}
\mor{a}{b}{$h$}
\enddc
$$

We have to prove that $(X,Cn_{\mathbb{P}}) \stackrel{\mathbb{P}(h)}{\longrightarrow} (X',Cn'_{\mathbb{P}})$ is, in fact, a morphism. That is, we have to verify that for all $A \subseteq X$, it holds that $h(Cn_{\mathbb{P}}(A)) = Cn'_{\mathbb{P}}(h(A))$, i.e., the following diagram commutes:

$$
\begindc{\commdiag}[120]
\obj(0,0)[c]{$\wp(X')$}
\obj(1,0)[d]{$\wp(X')$}
\obj(0,1)[a]{$\wp(X)$}
\obj(1,1)[b]{$\wp(X)$}
\mor{c}{d}{$Cn'_{\mathbb{P}}$}[-1,0]
\mor{a}{c}{$h$}[-1,0]
\mor{b}{d}{$h$}
\mor{a}{b}{$Cn_{\mathbb{P}}$}
\enddc
$$

The computation is straightforward, using lemma \ref{l21}:
$$
\begin{array}{lll}
h(Cn_{\mathbb{P}}(A)) & = & h(\bigcup \{Cn(A'): A' \subseteq A, Cn\text{-consistent}\}) \\
 & = & \bigcup \{h(Cn(A')): A' \subseteq A, Cn\text{-consistent}\} \\
 & = & \bigcup \{Cn'(h(A')): h(A') \subseteq h(A), Cn'\text{-consistent}\} \\
 & = & Cn'_{\mathbb{P}}(h(A)).
\end{array}
$$
The verification of functorial properties is immediate.
\end{prova}

\medskip

Let us examine which properties of the consequence operator $Cn$ are preserved by the functor $\mathbb{P}$.

\begin{definicao}
Let $Cn$ be a consequence operator on $X$ and $A,B \subseteq X$.

a) We say that $Cn$ satisfies \textbf{inclusion} iff (if and only if) $A \subseteq Cn(A)$;

b) We say that $Cn$ satisfies \textbf{idempotency} iff $Cn(Cn(A)) \subseteq Cn(A)$;

c) We say that $Cn$ satisfies \textbf{monotonicity} iff $A \subseteq B$ implies $Cn(A) \subseteq Cn(B)$;

d) We say that $Cn$ satisfies \textbf{finiteness} iff we have that:
$$
Cn(A) = \bigcup\{Cn(A'): A' \subseteq A \text{ finite}\}.
$$
Therefore, $x \in Cn(A)$ iff there exists a finite subset $A'$ of $A$ such that $x \in Cn(A')$.
\end{definicao}

\begin{proposicao}
\label{p32}
Let $(X,Cn)$ be a consequence structure and $A \subseteq X$. Then, it holds the following results:

a) If $A$ is $Cn$-consistent, then $Cn(A) \subseteq Cn_{\mathbb{P}}(A)$;

b) If $Cn$ is monotonic, then if $A$ is $Cn$-consistent, $Cn(A) = Cn_{\mathbb{P}}(A)$;
\end{proposicao}

\begin{prova}
a) Immediate from the definition of $Cn_{\mathbb{P}}$.

b) By part a) we have $Cn(A) \subseteq Cn_{\mathbb{P}}(A)$. Suppose that $x \in Cn_{\mathbb{P}}(A)$. Then, $x \in Cn(A')$ for some $Cn$-consistent subset $A'$ of $A$. By monotonicity, $Cn(A') \subseteq Cn(A)$. So, $x \in Cn(A)$.
\end{prova}

\begin{proposicao}
\label{p36}
The functor $\mathbb{P}$ enforces monotonicity. In other words, if $A \subseteq B$, then $Cn_{\mathbb{P}}(A) \subseteq Cn_{\mathbb{P}}(B)$.
\end{proposicao}

\begin{prova}
Suppose $A \subseteq B$ and $x \in Cn_{\mathbb{P}}(A)$. Then, there is some $Cn$-consistent subset $A' \subseteq A$ with $x \in Cn(A)$. But that same $A'$ is a $Cn$-consistent subset of $B$ with $x \in Cn(A)$. So, $x \in Cn_{\mathbb{P}}(B)$.
\end{prova}

\begin{proposicao}
\label{p35}
The functor $\mathbb{P}$ preserves finiteness, i.e., if $Cn$ satisfies finiteness, then $Cn_{\mathbb{P}}$ also satisfies finiteness.
\end{proposicao}

\begin{prova}
Suppose that $Cn$ satisfies finiteness. Let $x \in Cn_{\mathbb{P}}(A)$. By definition of $Cn_{\mathbb{P}}$, there exists $A' \subseteq A$, $Cn$-consistent, such that $x \in Cn(A')$. So, there is a finite $A^{*} \subseteq A'$ such that $x \in Cn(A^{*})$. Since $Cn$ satisfies finiteness, $Cn(A^{*} \subseteq Cn(A"')$, so $A^{*}$ is $Cn$-consistent as well. By proposition \ref{p32} (a), since $A^{*}$ is $Cn$-consistent, we have $Cn(A^{*}) \subseteq Cn_{\mathbb{P}}(A^{*})$, so $x \in Cn_{\mathbb{P}}(A^{*})$. Therefore, there exists $A^{*} \subseteq A$, finite, such that $x \in Cn_{\mathbb{P}}(A^{*})$.

On the other hand, suppose that there is a subset $A'$ of $A$, finite, such that $x \in Cn_{\mathbb{P}}(A')$. 
So, by proposition \ref{p36}, $Cn_{\mathbb{P}}(A') \subseteq Cn_{\mathbb{P}}(A)$, and $x \in Cn_{\mathbb{P}}(A)$.
\end{prova}

\medskip
We say that a consequence structure $(X,Cn)$ is \emph{normal} (or \emph{Tarskian}) if and only if the consequence operator satisfies inclusion, idempotency and monotonicity.

\begin{proposicao}
\label{p37}
If $(X,Cn)$ is normal and there is $u \in X$ such that $\{u\}$ is $Cn$-inconsistent, then in $(X,Cn_{\mathbb{P}})$ there is no $Cn_{\mathbb{P}}$-inconsistent sets.
\end{proposicao}

\begin{prova}
Consider the hypotheses and suppose, \emph{ad absurdum}, $A \subseteq X$ such that $Cn_{\mathbb{P}}(A) = X$. So, $u \in Cn_{\mathbb{P}}(A)$. Therefore, there is $A' \subseteq A$, $Cn$-consistent such that $u \in Cn(A')$. Hence, $\{u\} \subseteq Cn(A')$ and we have $X = Cn(\{u\}) \subseteq Cn(Cn(A')) = Cn(A')$. (contradiction!)
\end{prova}

\begin{corolario}
In the conditions of the proposition above, we have 
$$
\mathbb{P}(\mathbb{P}(X,Cn)) = \mathbb{P}(X,Cn) = (X,Cn_{\mathbb{P}})
$$
i.e., the functor acts in an idempotent way.
\end{corolario}

\begin{prova}
By proposition \ref{p37}, every subset $A$ of $X$ is $Cn_{\mathbb{P}}$-consistent. Therefore, by proposition \ref{p32} (b), $Cn_{\mathbb{P}}(Cn_{\mathbb{P}}(A)) = Cn_{\mathbb{P}}(A)$.
\end{prova}

\medskip
We will see, in the section 5, that the functor $\mathbb{P}$ preserves neither inclusion nor idempotency.\footnote{If we want to enforce inclusion, we could make a slight modification in the definition of $\mathbb{P}$ in order to include $A$ into $Cn_{\mathbb{P}}(A)$, such that
$$
Cn_{\mathbb{P}}(A) := A \cup \bigcup\{Cn(A'): A' \subseteq A,\   Cn\text{-consistent} \}.
$$
}


\section{Paraconsistent consequence structures}
We will show that a paraconsistent transformation indeed turns structures that satisfy \emph{ex falso quodlibet} into a consequence structure in which this principle fails. In order to do so, we have to introduce some conventions and definitions.\footnote{Chakraborty and Dutta in \cite{chakraborty} explore a way to axiomatize paraconsistent consequence structures
using abstract consequence operators.}

From now on, we suppose that the set $X$ is endowed with an operator intended to be a \emph{negation operation}, denoted by the symbol $\neg$. Thus, if $x \in X$, then $\neg x \in X$ and $\neg x$ is called the \emph{negation} of $x$.

\begin{definicao}
Let $(X,Cn)$ be a consequence structure.

1. We say that $(X,Cn)$ satisfies \textbf{ex falso quodlibet} (or satisfies \textbf{explosion}, or is \textbf{explosive}) iff for all $A \subseteq X$, if there is $x \in X$ such that $x, \neg x \in Cn(A)$, then $Cn(A) = X$ (i.e., $A$ is $Cn$-inconsistent). Otherwise, $(X,Cn)$ is called \textbf{paraconsistent}.

2. We say that $(X,Cn)$ satisfies \textbf{joint consistency} iff there exists $x \in X$ such that $\{ x \},\{\neg x \}$ are both $Cn$-consistent and $\{x,\neg x\}$ is $Cn$-inconsistent.

3. We say that $(X,Cn)$ satisfies \textbf{conjunctive property} iff for all $x,y \in X$, there exists $z \in X$ such that $Cn(\{ x,y \}) = Cn(\{ z \})$.
\end{definicao}

Now, we present a sufficient condition for the paraconsistentization functor $\mathbb{P}$  to transform a consequence structure into a paraconsistent one. Therefore, that functor deserves its name!

\begin{teorema}
If $(X,Cn)$ is normal, explosive, satisfies joint consistency and also conjunctive property, then $(X,Cn_{\mathbb{P}})$ is paraconsistent.
\end{teorema}

\begin{prova}
Since $(X,Cn)$ satisfies joint consistency, there exists $a \in X$ such that $\{ a \}, \{ \neg a \}$ are both $Cn$-consistent. Consider $A = \{ a, \neg a \}$ and, then, $Cn(A) = X$. By inclusion, $A \subseteq Cn(A)$ and by joint consistency, $A \subseteq Cn_{\mathbb{P}}(A)$. As $(X,Cn)$ satisfies conjunctive property, there exists $c \in X, Cn(\{ c \}) = Cn(A) = X$. We will show that $c \notin Cn_{\mathbb{P}}(A)$, i.e., $(X,Cn_{\mathbb{P}})$ is paraconsistent. The set $A$ has three $Cn$-consistent subsets, that is: $\{ a \}$, $\{ \neg a \}$ and $\emptyset$ (the empty set). ($\emptyset$ is $Cn$-consistent for $\emptyset \subseteq \{ a \}$ and, by monotonicity, $Cn(\emptyset) \subseteq Cn(\{ a \}) \neq X$.) We have to show that $c$ does not belong to the operator $Cn$ applied to these sets.
If $c \in Cn(\{ a \})$, then $\{ c \} \subseteq Cn(\{ a \})$. By monotonicity,  $Cn(\{ c \}) \subseteq Cn(Cn(\{ a \}))$. By idempotency and inclusion, $X \subseteq Cn(\{ c \}) \subseteq Cn(\{ a \})$, i.e., $Cn(\{ a \}) = X$ (contradiction!). The same argument shows that $c \notin Cn(\{ \neg a \})$. Finally, if $c \in Cn(\emptyset)$, then $\{ c \} \subseteq Cn(\emptyset)$. So, we have, $X = Cn(\{ c \}) \subseteq Cn(Cn(\emptyset)) = Cn(\emptyset)$ (contradiction!); therefore, $(X,Cn_{\mathbb{P}})$ is paraconsistent. 
\end{prova}

\medskip
In the next section, we study a particular case of paraconsistentization.


\section{Paraconsistentization of propositional classical logic}

Let $X$ be the set of formulas of an usual propositional language with $\neg$ (negation), $\vee$ (disjunction), $\wedge$ (conjunction), $\rightarrow$ (implication) and propositional letters: $p,q,r,...,p_{1},q_{1},r_{1},...$ and so on. Let $Cn$ be the standard consequence operator of the propositional classical logic.

Let $(X,Cn_{\mathbb{P}})$ be the result of the action of $\mathbb{P}$ on $(X,Cn)$. We call $(X,Cn_{\mathbb{P}})$ a \emph{propositional paraclassical logic}. We proceed to study some properties of $(X,Cn_{\mathbb{P}})$.\footnote{Details of the following constructions can be found in \cite{desouza}, but not using the framework
provided by the paraconsistentization functor introduced here.} 

For convenience, we use $A \vdash a$ to denote $a \in Cn(A)$. On the contrary, we use $A \nvdash a$. Similarly for $A \vdash_{\mathbb{P}} a$.

\begin{proposicao}
Inclusion and idempotency do not hold in $(X,Cn_{\mathbb{P}})$. Therefore, these properties are not preserved by the functor $\mathbb{P}$.
\end{proposicao}

\begin{prova}
For inclusion, notice that $\{p \wedge \neg p \} \nvdash_{\mathbb{P}} p \wedge \neg p$ because $\emptyset$ is the only $Cn$-consistent subset of $\{p \wedge \neg p \}$ and $\emptyset \nvdash p \wedge \neg p$.	

For idempotency, let $A = \{p, \neg p\}$. Then, $p \vee q , \neg p \in Cn_{\mathbb{P}}(A)$. Therefore, $q \in Cn_{\mathbb{P}}(Cn_{\mathbb{P}}(A))$, but $q \notin Cn_{\mathbb{P}}(A)$.
\end{prova}

\medskip
It is well known that $(X,Cn)$ satisfies the property of \emph{transitivity}: if $A \vdash b$ for every $b \in B$ and $B \vdash a$, then $A \vdash a$.
Moreover, we have a \emph{weak form of transitivity}: if $A \vdash b$ and $\{b\} \vdash c$, then $A \vdash c$.

\begin{proposicao}
Transitivity does not hold in $(X,Cn_{\mathbb{P}})$.
\end{proposicao}

\begin{prova}
Consider $A = \{p, \neg p \}$, $B = \{ p \vee q, \neg p \}$ and $a = q$. Then, $\{p, \neg p \} \vdash_{\mathbb{P}} p \vee q$, $\{p, \neg p \} \vdash_{\mathbb{P}} \neg p$, $\{p \vee q, \neg p \} \vdash_{\mathbb{P}} q$, but $\{p, \neg p \} \nvdash_{\mathbb{P}}q$.
\end{prova}

\medskip
For weak transitivity, we need a preliminary result. In $(X,Cn)$, we say that $a \in X$ is a $Cn$-\emph{contradiction} iff $\{ a  \} \vdash p \wedge \neg p$. Moreover, $a$ is a $Cn$-\emph{theorem} iff $\emptyset \vdash a$.

\begin{lema}
\label{l53}
We have the following properties:

a) If $b$ is a $Cn$-contradiction, then for every $A \subseteq X$, $A \nvdash_{\mathbb{P}} b$;

b) If $b$ is a $Cn$-theorem and $\{b\} \vdash_{\mathbb{P}} c$, then $c$ is a $Cn$-theorem and for every $A \subseteq X$, $A \vdash_{\mathbb{P}} c$;

c) If $\{a\} \vdash_{\mathbb{P}} b$, then $b$ is a $Cn$-theorem or $\{a\}$ is $Cn$-consistent and $\{a\} \vdash b$.
\end{lema}

\begin{prova}
Immediate from the definitions.
\end{prova}

\begin{proposicao}
Weak transitivity holds in $(X,Cn_{\mathbb{P}})$. That is, if $A \vdash_{\mathbb{P}} b$ and $\{b\} \vdash_{\mathbb{P}} c$, then $A \vdash_{\mathbb{P}} c$.
\end{proposicao}

\begin{prova}
Since $\{b\} \vdash_{\mathbb{P}} c$, there are some $Cn$-consistent $K \subseteq \{b\}$ with $K \vdash c$. Either $K = \emptyset$ or $K = \{b\}$; these are the only two subsets of $\{b\}$. In the first case, $c$ is a $Cn$-theorem, and so $A \vdash c$, by the $Cn$-consistency of $\emptyset$ and monotonicity of $Cn_{\mathbb{P}}$. In the second case, $\{b\} \vdash c$. Since $A \vdash_{\mathbb{P}} b$, there are some $Cn$-consistent $A' \subseteq A$ with $A' \vdash b$. By the transitivity of $Cn$, this gives $A' \vdash c$. So, there is some $Cn$-consistent $A' \subseteq A$ with $A' \vdash c$; that is, $A \vdash_{\mathbb{P}} c$.\footnote{We thank one more time one of the referees for suggesting this corrected version of the proof.}
\end{prova}

\medskip
In $(X,Cn)$, it holds the \emph{deduction theorem}\footnote{Notice that we are considering the deduction theorem in one direction only.}: if $A \cup \{a\} \vdash b$, then $A \vdash a \rightarrow b$.

\begin{proposicao}
\label{p55}
The deduction theorem is valid in $(X,Cn_{\mathbb{P}})$.
\end{proposicao}

\begin{prova}
Suppose that $A \cup \{a\} \vdash_{\mathbb{P}} b$. Then, there exists $B \subseteq A \cup \{a\}$, $Cn$-consistent, such that $B  \vdash b$. We have two cases:

1) $B \subseteq A$. In this case,  $B \cup \{a\} \vdash b$ and $B \vdash a \rightarrow b$, by deduction theorem for $(X,Cn)$. Since $B$ is $Cn$-consistent, we have that $A \vdash_{\mathbb{P}} a \rightarrow b$.

2) $B \nsubseteq A$. In this case, we have $B - \{a\} \vdash a \rightarrow b$ and $B - \{a\}$ is $Cn$-consistent. So, $A \vdash_{\mathbb{P}} a \rightarrow b$.

And this completes the proof.
\end{prova}

\medskip
Notice that \emph{modus ponens} does not hold in $(X,Cn_{\mathbb{P}})$, for example in the case that the conclusion of the rule would be a $Cn$-contradiction. By the same reason, the converse of deduction theorem does not hold in $(X,Cn_{\mathbb{P}})$. For example, we have $A \vdash_{\mathbb{P}} (p \wedge \neg p) \rightarrow (p \wedge \neg p)$, but $A \cup \{p \wedge \neg p \} \nvdash_{\mathbb{P}} p \wedge \neg p$.

A set $A \subseteq X$ is called $Cn$-\emph{contradictory} iff there is a formula $a$ such that $A \vdash a$ and $A \vdash \neg a$. We say that $A$ is $Cn$-\emph{strongly contradictory} iff there is a $Cn$-contradictory formula $a$ such that $A \vdash a$. Moreover, we say that $A$ is $Cn$-\emph{paraconsistent} iff $A$ is $Cn$-consistent and $Cn$-contradictory. The same definitions hold for $Cn_{\mathbb{P}}$.

In $(X,Cn)$, there are no $Cn$-paraconsistent sets. On the other hand, in $(X,Cn_{\mathbb{P}})$, we have $Cn_{\mathbb{P}}$-paraconsistent sets but we do not have neither $Cn_{\mathbb{P}}$-inconsistent sets nor $Cn_{\mathbb{P}}$-strongly contradictory sets.

We can summarize the results of this section in the following table.

\begin{center}
\begin{tabular}{|c|c|c|}
\hline
 & $Cn$ & $Cn_{\mathbb{P}}$ \\
\hline
finiteness & \checkmark & \checkmark \\
monotonicity & \checkmark & \checkmark \\
inclusion & \checkmark  & $\times$ \\
idempotency &  \checkmark & $\times$  \\
transitivity & \checkmark  & $\times$ \\
weak transitivity & \checkmark  & \checkmark \\
deduction & \checkmark  & \checkmark \\
inconsistent sets & \checkmark  & $\times$  \\
contradictory sets & \checkmark  & \checkmark  \\
strongly contradictory sets & \checkmark  & $\times$  \\
paraconsistent sets & $\times$  & \checkmark  \\
\hline
\end{tabular}

\medskip
\checkmark means that the operator of the consequence structure has the property or there are sets as indicated.

$\times$ means the contrary.
\end{center}


\section{Conclusion}

Newton da Costa conjectured in many different places that all logics can be adapted in order to become
paraconsistent. The device used in the present work realizes this task. All theories, regardless of their nature, require 
an underlying logic. In most cases,
this logic is classical and, thus, contradictions are not allowed. The price to pay for
finding a contradiction is sometimes too high that the theory has to be abandoned or reformulated. 
Notwithstanding, paraconsistentizing the underlying logic can save a given
theory. Basically, methods of paraconsistentization have applications wherever
paraconsistency plays a role: solving epistemic paradoxes, dealing with deontic dilemmas,
modelling inconsistent reasoning in general, and everywhere we need 
logics for underlying contradictory but non-explosive theories, we can apply the methodology developed here.
Therefore, the methodology proposed has a very large range of application in contexts involving contradictions.

We have developed a way to convert a given logic into a paraconsistent one. In particular,
we focused on doing this by way of functors defined in categories where consequence structures 
are objects. This kind of approach can be featured 
in universal logic. Instead of exploring particular
logical systems, universal logic investigates all possible logics (see \cite{beziau2}). 
As pointed out in \cite{costaleite}, 
there is no unique method for paraconsistentizing
a given logic, but rather a plurality of them.
In this sense,
paraconsistentization is for paraconsistent logics what universal logic is for logics
in general. Therefore, it is a general theory of paraconsistent logics. 
Indeed, there are many methods and ways one can use in order to paraconsistentize a given
non-explosive logic. These other possibilities should still be studied.

Our approach
to paraconsistency does not coincide with other researches especially because our
paraconsistentization is realized without explicit mention to the concept
of \emph{negation}. It produces an unexpected result: \emph{modus ponens} is not
generally valid.\footnote{Ja\'skowski in \cite{jaskowski} shows, for example, a case
in which \emph{modus ponens} fails: classical logic is not able to be an
underlying logic of a discussive system.} This can be a contribution to discussions
regarding the nature of paraconsistency and its relation with negation.\footnote{Negation
is required, however, to show that the functor does what it is supposed to do. But this
is a different issue.} In addition,
it offers another way to deal with inconsistency. Suppose a logician trying to formalize
a given theory and investigating its logical consequences. Then, an explosive logic - which
seems to be adequate for some tasks - is chosen. Nonetheless, during the investigation,
contradictions are found and they cause trivialization of the system. By applying
a paraconsistentization functor, it is possible to restrict the domain of the original explosive
logic to the consistent subsets of the theory. The advantage os this approach is
that it is able to keep main characteristics of the input logic and, in particular, its
theorems. Thus, for example, an intuitionist logician can still study a constructive theory -
and yet inconsistent - using the paraconsistent counterpart of it, hence keeping crucial
features of the intuitionistic logic but without trivialization.

Last, but not least, the work developed by Rescher and Manor (in \cite{rescher-manor}) presents logical 
machinery for obtaining non-trivial consequences from inconsistent sets. Our task is not a 
generalization of their work. There are essential differences. Notice that the notion of consistency in 
our paper is independent of any underlying language and, in particular, it is independent of the concept of 
negation, while their characterization of consistency relies on a previously established language.
In addition, their paper contains a very strong presupposition: it is not rational to accept all 
consequences of inconsistent premises (p.182). 
This reasoning presupposes exactly 
one of the core philosophical tenets that paraconsistent logics try to overcome, namely: the
parochial thesis according to which  consistency 
is a necessary condition for rationality (for a survey of the philosophical implications of paraconsistency, 
see \cite{dias}). Once a paraconsistent logic is formulated, the consequences of inconsistent 
premises are no longer trivial, so they can all be accepted rationally.

In conclusion, future lines of research could include exploring other universal properties holding in
the category $CON$, as well as examining more properties preserved (or lost) when the
paraconsistentization functor is applicable. In a less abstract mode, we can think about consequence
structures endowed with syntactical and semantical dimensions; and, at this level,
questions of whether metalogical properties like soundness and completeness are preserved
by paraconsistentization could also be examined.

\subsection*{Acknowledgement}
Thanks to anonymous reviewers for a detailed analysis of
the paper which helped us to improve it. Thanks also to
Scott Randall Paine for proof-reading the English.

\end{document}